\documentclass[11pt]{article}
\usepackage{a4wide}

\newcommand{\ba}{\begin{array}}
\newcommand{\ea}{\end{array}}
\newcommand{\be}{\begin{equation}}
\newcommand{\ee}{\end{equation}}
\newcommand{\la}{\label}
\newcommand{\bea}{\begin{eqnarray}}
\newcommand{\eea}{\end{eqnarray}}
\newcommand{\ch}{\choose}
\renewcommand{\a}{\alpha}
\renewcommand{\b}{\beta}

\newcommand{\G}{\Gamma}

\newcommand{\KL}{L_n^{\a,M}(x)}

\newcommand{\GP}{P_n^{\a,\b,M,N}(x)}
\newcommand{\SP}{P_n^{(\a,\a)}(x)}
\newcommand{\SGP}{P_n^{\a,\a,M,M}(x)}
\renewcommand{\l}{\left}
\renewcommand{\r}{\right}
\newcommand{\set}[1]{\left\{#1\right\}_{n=0}^{\infty}}
\newcommand{\hyp}[5]{\mbox{}_{#1}F_{#2}
\left(\left.\begin{array}{c}#3\\#4\end{array}\right|#5\right)}
\newcommand{\ndots}{n=0,1,2,\ldots}
\newcommand{\n}{\nonumber}
\newcommand{\nn}{\nonumber \\}
\newcommand{\ds}{\displaystyle}
\newcounter{stelling}
\newcommand{\st}[1]{\par\vspace{0.5cm}\refstepcounter{stelling}
{\bf Theorem \thestelling.} {\sl #1}\par\vspace{0.5cm}}

\begin{document}

\begin{center}
{\Large DIFFERENTIAL EQUATIONS\\ FOR\\ SYMMETRIC GENERALIZED\\
\vspace{2mm}ULTRASPHERICAL POLYNOMIALS}

\vspace{1cm}

{\large Roelof Koekoek}

\vspace{1cm}

Delft University of Technology\\ Faculty of Technical Mathematics and Informatics\\
Mekelweg 4\\ 2628 CD Delft\\ The Netherlands\\ e-mail : koekoek@twi.tudelft.nl
\end{center}

\vspace{1cm}

\begin{abstract}
We look for differential equations satisfied by the generalized Jacobi
polynomials $\set{\GP}$ which are orthogonal on the interval $[-1,1]$ with
respect to the weight function
$$\frac{\G(\a+\b+2)}{2^{\a+\b+1}\G(\a+1)\G(\b+1)}(1-x)^{\a}(1+x)^{\b}+
M\delta(x+1)+N\delta(x-1),$$
where $\a>-1$, $\b>-1$, $M\ge 0$ and $N\ge 0$.

In the special case that $\b=\a$ and $N=M$ we find all differential equations
of the form
$$\sum_{i=0}^{\infty}c_i(x)y^{(i)}(x)=0,\;y(x)=\SGP,$$
where the coefficients $\l\{c_i(x)\r\}_{i=1}^{\infty}$ are independent of
the degree $n$.

We show that if $M>0$ only for nonnegative integer values of $\a$ there
exists exactly one differential equation which is of finite order $2\a+4$.

By using quadratic transformations we also obtain differential equations for
the polynomials $\set{P_n^{\a,\pm\frac{1}{2},0,N}(x)}$ for all $\a>-1$ and
$N\ge 0$.
\end{abstract}

\vspace{1cm}

{\bf AMS Subject Classification :} Primary 33C45 (33A65), Secondary 34A35

\newpage

\section{Introduction}

In the late thirties (see \cite{Krall38} and \cite{Krall40}) H.L. Krall
classified all sets of orthogonal polynomials $\set{P_n(x)}$ with
degree$[P_n(x)]=n$ which satisfy a fourth order differential equation of
the form
$$p_4(x)y^{(4)}(x)+p_3(x)y^{(3)}(x)+p_2(x)y''(x)+p_1(x)y'(x)+p_0(x)y(x)=0$$
where $\l\{p_i(x)\r\}_{i=0}^{4}$ are polynomials with degree$[p_i(x)]\le i$
and $\l\{p_i(x)\r\}_{i=1}^4$ are independent of the degree $n$.
These sets of orthogonal polynomials include the classical Legendre, Laguerre,
Hermite, Bessel and Jacobi polynomials. He also found three other sets of
orthogonal polynomials satisfying a fourth order differential equation of
this type. In \cite{Krall81} A.M. Krall studied these new sets of orthogonal
polynomials in more details and named them the Legendre type, Laguerre type
and Jacobi type polynomials.
These polynomials are generalizations of the classical Legendre, Laguerre
(with $\a=0$) and Jacobi polynomials (with $\b=0$) in the sense that the
weight function for these orthogonal polynomials consists of the classical
weight function together with a Dirac delta function at the endpoint(s) of
the interval of orthogonality.

Later L.L. Littlejohn (see \cite{Lit_Krall}) studied a generalization of
the Legendre type polynomials and named them after H.L. Krall~:
the Krall polynomials. These Krall polynomials are orthogonal on the interval
$[-1,1]$ with respect to the weight function
$$\frac{1}{A}\delta(x+1)+\frac{1}{B}\delta(x-1)+C,
\;A>0,\;B>0\;\mbox{ and }\;C>0.$$
In general ($A\ne B$), these polynomials do not fit in the class of
polynomials which satisfy a fourth
order differential equation of the above type. The Krall polynomials satisfy a
sixth order differential equation of a similar form.

A.M. Krall and L.L. Littlejohn did some work on the classification of higher
order differential equations having orthogonal polynomial solutions. They
tried to classify all differential equations of the form
$$\sum_{i=0}^rp_i(x)y^{(i)}(x)=0,\;r\in\{2,3,4,\ldots\},$$
where $\l\{p_i(x)\r\}_{i=0}^r$ are polynomials with degree$[p_i(x)]\le i$ and
$\l\{p_i(x)\r\}_{i=1}^r$ are independent of $n$ having orthogonal polynomial
solutions $\set{P_n(x)}$ with degree$[P_n(x)]=n$.
See \cite{ClassI} and \cite{ClassII}.

In \cite{Koorn} T.H. Koornwinder found a general class of orthogonal polynomials
which generalize the Legendre type, Jacobi type and Krall polynomials. These
polynomials are orthogonal on the interval $[-1,1]$ with respect to the weight
function
$$\frac{\G(\a+\b+2)}{2^{\a+\b+1}\G(\a+1)\G(\b+1)}
(1-x)^{\a}(1+x)^{\b}+M\delta(x+1)+N\delta(x-1),$$
where $\a>-1$, $\b>-1$, $M\ge 0$ and $N\ge 0$. For these generalized Jacobi
polynomials we will use Koornwinder's notation~: $\set{\GP}$.

As a limit case he found the polynomials $\set{\KL}$ which are orthogonal
on the interval $[0,\infty)$ with respect to the weight function
$$\frac{1}{\G(\a+1)}x^{\a}e^{-x}+M\delta(x),\;\a>-1\;\mbox{ and }\;M\ge 0.$$
These polynomials generalize the classical Laguerre polynomials.

In \cite{DV} J. Koekoek and R. Koekoek showed that the polynomials $\set{\KL}$
satisfy a unique differential equation of the form
$$M\sum\limits_{i=0}^{\infty}a_i(x)y^{(i)}(x)+xy''(x)+(\a+1-x)y'(x)+ny(x)=0,$$
where $\l\{a_i(x)\r\}_{i=0}^{\infty}$ are continuous functions on the real
line and $\l\{a_i(x)\r\}_{i=1}^{\infty}$ are independent of $n$.
It turns out that the coefficients $\l\{a_i(x)\r\}_{i=0}^{\infty}$ are polynomials
and the differential equation is of infinite order in general if $M>0$. However,
only for nonnegative integer values of $\a$ the order reduces to $2\a+4$.

We note that it is well-known that all sets of polynomials named before satisfy
a second order differential equation with polynomial coefficients depending on $n$, but
of bounded degree. See for instance \cite{Koorn}, \cite{Lit_Krall2dv} and
\cite{Lit&Shore}.

In this paper we look for differential equations for the polynomials
$\set{\GP}$ with $\b=\a$ and $N=M$. Until now, only two special
cases, due to H.L.~Krall and L.L.~Littlejohn, are known.
In \cite{Krall38} H.L.~Krall showed that the polynomials $\set{P_n^{0,0,M,M}(x)}$
satisfy the following fourth order (if $M>0$) differential equation~:
\bea & &-\frac{1}{2}M(1-x^2)^2y^{(4)}(x)+4Mx(1-x^2)y^{(3)}(x)+{}\nn
& &\hspace{1cm}{}+(6M+1)(1-x^2)y''(x)-2xy'(x)
+\frac{1}{2}n(n+1)\l[(n-1)(n+2)M+2\r]y(x)=0\n\eea
and later L.L.~Littlejohn found the following sixth order (if $M>0$)
differential equation for the polynomials $\set{P_n^{1,1,M,M}(x)}$~:
\bea & &\frac{1}{36}M(1-x^2)^3y^{(6)}(x)-\frac{2}{3}Mx(1-x^2)^2y^{(5)}(x)
-\frac{5}{3}M(1-x^2)(1-3x^2)y^{(4)}(x)+{}\nn
& &\hspace{1cm}{}+\frac{40}{3}Mx(1-x^2)y^{(3)}(x)+(10M+1)(1-x^2)y''(x)+{}\nn
& &\hspace{1cm}{}-4xy'(x)+\frac{1}{36}n(n+3)
\l[(n-1)(n+1)(n+2)(n+4)M+36\r]y(x)=0\n\eea
both in a different notation. The latter sixth order differential equation
has not appeared in the literature yet.

We will derive all differential equations for the polynomials $\set{\SGP}$
for every $\a>-1$ and $M\ge 0$, which are of the form
$$\sum_{i=0}^{\infty}c_i(x)y^{(i)}(x)=0,$$
where the coefficients $\l\{c_i(x)\r\}_{i=0}^{\infty}$ are continuous
functions on the real line and $\l\{c_i(x)\r\}_{i=1}^{\infty}$ are
independent of the degree $n$.

So, we consider the polynomials $\set{\SGP}$ which can be defined by
(see \cite{Koorn}, in a slightly different notation)
\be\la{def}\SGP=C_0\SP-C_1x\frac{d}{dx}\SP,\ee
where
\be\la{C}\l\{\ba{l}\ds C_0=1+M\frac{2n}{(\a+1)}{n+2\a+1 \ch n}+
4M^2{n+2\a+1 \ch n-1}^2\\ \\
\ds C_1=\frac{2M}{(2\a+1)}{n+2\a \ch n}+
\frac{2M^2}{(\a+1)}{n+2\a \ch n-1}{n+2\a+1 \ch n}.\ea\r.\ee
As Koornwinder already remarked (see \cite{Koorn}) the case $2\a+1=0$ must
be understood by continuity in $\a$.

Further we will show that for $M>0$ these differential equations are of
infinite order in general and only for nonnegative integer values of $\a$ we
find exactly one differential equation of finite order $2\a+4$. This answers one
of the questions raised in \cite{Conj} by W.N.~Everitt and L.L.~Littlejohn.

Finally, we will also derive differential equations for the polynomials
$\set{P_n^{\a,\pm\frac{1}{2},0,N}(x)}$ for all $\a>-1$ and $N\ge 0$.

\section{Some classical formulas}

In this section we give the definition and some properties of the classical
ultraspherical polynomials $\set{\SP}$. For details the reader is referred
to \cite{Chi} and \cite{Szego}. We will only give those properties we need in
this paper. Further we list some other classical formulas which we will use
later on.

The polynomials $\set{\SP}$ can be defined by their representation as a
hypergeometric function as
\be\la{defJ}\SP={n+\a \ch n}\hyp{2}{1}{-n,n+2\a+1}{\a+1}{\frac{1-x}{2}},\;\ndots.\ee

A simple consequence of this definition is the differentiation formula
\bea\la{afgJ}& &D^i\SP\nn
&=&{n+\a \ch n}\l(-\frac{1}{2}\r)^i\sum_{k=i}^{\infty}
\frac{(-n)_k(n+2\a+1)_k}{(k-i)!(\a+1)_k}\l(\frac{1-x}{2}\r)^{k-i}\nn
&=&{n+\a \ch n}\l(-\frac{1}{2}\r)^i\sum_{k=0}^{\infty}
\frac{(-n)_{k+i}(n+2\a+1)_{k+i}}{k!(\a+1)_{k+i}}\l(\frac{1-x}{2}\r)^k,
\;i=0,1,2,\ldots.\eea

Further we have the well-known symmetry relation
\be\la{symJ}P_n^{(\a,\a)}(-x)=(-1)^n\SP,\;\ndots.\ee

The ultraspherical polynomials satisfy a second order linear differential
equation given by
$$(1-x^2)y''(x)-2(\a+1)xy'(x)+n(n+2\a+1)y(x)=0.$$
By using induction it is easy to show that this differential equation
implies that
\bea\la{dvJ}& &(1-x^2)D^{i+2}\SP-2(\a+i+1)xD^{i+1}\SP+{}\nn
& &\hspace{3cm}{}+(n-i)(n+2\a+i+1)D^i\SP=0,\;i=0,1,2,\ldots.\eea

We also need the following formulas for the even and odd order
ultraspherical polynomials~:
\be\la{even}P_{2n}^{(\a,\a)}(x)=\l(-\frac{1}{4}\r)^n{2n+\a \ch n}
\hyp{2}{1}{-n,n+\a+\frac{1}{2}}{\frac{1}{2}}{x^2},\;\ndots\ee
and
\be\la{odd}P_{2n+1}^{(\a,\a)}(x)=\l(-\frac{1}{4}\r)^n{2n+\a \ch n}
(2n+\a+1)x\hyp{2}{1}{-n,n+\a+\frac{3}{2}}{\frac{3}{2}}{x^2},\;\ndots,\ee
respectively. These formulas can be found in \cite{Bate} in a slightly
different notation.

We will often use the Vandermonde summation formula~:
\be\la{Van}\hyp{2}{1}{-n,b}{c}{1}=\frac{(c-b)_n}{(c)_n},\;\ndots\ee
and the Saalsch\"{u}tz summation formula~:
\be\la{Saal}\hyp{3}{2}{-n,a,b}{c,-n+a+b-c+1}{1}=
\frac{(c-a)_n(c-b)_n}{(c)_n(c-a-b)_n},\;\ndots.\ee
These summation formulas can be found in \cite{Bate} for instance.

Finally, we remark that the Taylor series at the point zero of a
hypergeometric function of the form $_{p+1}F_p$ has a radius of
convergence $1$ unless it terminates.
Moreover, such a series also converges absolutely at $1$ if the sum of the
numerator parameters is less than the sum of the denominator parameters.
For details the reader is referred to \cite{Brom}.

\section{The differential equations for $\SGP$}

In order to find all differential equations of the form
\be\la{dv}M\sum_{i=0}^{\infty}a_i(x)y^{(i)}(x)+(1-x^2)y''(x)-2(\a+1)xy'(x)
+n(n+2\a+1)y(x)=0\ee
for the polynomials $\set{\SGP}$, where the coefficients
$\l\{a_i(x)\r\}_{i=0}^{\infty}$ are continuous functions on the real line
and where $\l\{a_i(x)\r\}_{i=1}^{\infty}$ are independent of $n$, we set
$y(x)=\SGP$ in (\ref{dv}) and use the definition (\ref{def}) to find
\bea & &MC_0\sum_{i=0}^{\infty}a_i(x)D^i\SP-MC_1\sum_{i=0}^{\infty}ia_i(x)D^i\SP+{}\nn
& &\hspace{2cm}{}-MC_1x\sum_{i=0}^{\infty}a_i(x)D^{i+1}\SP+{}\nn
& &{}+(1-x^2)\l[C_0\frac{d^2}{dx^2}\SP-2C_1\frac{d^2}{dx^2}\SP-C_1x\frac{d^3}{dx^3}\SP\r]+{}\nn
& &{}-2(\a+1)x\l[C_0\frac{d}{dx}\SP-C_1\frac{d}{dx}\SP-C_1x\frac{d^2}{dx^2}\SP\r]+{}\nn
& &{}+n(n+2\a+1)\l[C_0\SP-C_1x\frac{d}{dx}\SP\r]=0.\n\eea
Now we use (\ref{dvJ}) with $i=0$ and $i=1$ to obtain
\bea & &MC_0\sum_{i=0}^{\infty}a_i(x)D^i\SP-MC_1\l[\sum_{i=0}^{\infty}ia_i(x)D^i\SP
+x\sum_{i=0}^{\infty}a_i(x)D^{i+1}\SP\r]\nn
&=&2C_1\frac{d^2}{dx^2}\SP.\n\eea
We consider both sides as polynomials in $M$. Comparing
the coefficients of equal powers of $M$ on both sides leads to
\bea &M:& \sum_{i=0}^{\infty}a_i(x)D^i\SP=\frac{4}{(2\a+1)}{n+2\a \ch n}
\frac{d^2}{dx^2}\SP\nn
&M^2:& \frac{2n}{(\a+1)}{n+2\a+1 \ch n}\sum_{i=0}^{\infty}a_i(x)D^i\SP+{}\nn
& &\hspace{1cm}{}-\frac{2}{(2\a+1)}
{n+2\a \ch n}\l[\sum_{i=0}^{\infty}ia_i(x)D^i\SP+x\sum_{i=0}^{\infty}a_i(x)D^{i+1}\SP\r]\nn
& &=\frac{4}{(\a+1)}{n+2\a \ch n-1}{n+2\a+1 \ch n}\frac{d^2}{dx^2}\SP\nn
&M^3:& 4{n+2\a+1 \ch n-1}^2\sum_{i=0}^{\infty}a_i(x)D^i\SP\nn
& &=\frac{2}{(\a+1)}{n+2\a \ch n-1}{n+2\a+1 \ch n}\times{}\nn
& &\hspace{2cm}{}\times\l[\sum_{i=0}^{\infty}ia_i(x)D^i\SP
+x\sum_{i=0}^{\infty}a_i(x)D^{i+1}\SP\r].\n\eea
It is clear that the third relation is a linear combination of the first two
relations which can be simplified to
\be\la{rel1}\sum_{i=0}^{\infty}a_i(x)D^i\SP=\frac{4}{(2\a+1)}{n+2\a \ch n}
\frac{d^2}{dx^2}\SP\ee
and
\be\la{rel2}\sum_{i=0}^{\infty}ia_i(x)D^i\SP+x\sum_{i=0}^{\infty}a_i(x)D^{i+1}\SP
=4{n+2\a+1 \ch n-1}\frac{d^2}{dx^2}\SP.\ee

Since we demand that the coefficients $\l\{a_i(x)\r\}_{i=1}^{\infty}$ are
independent of $n$ we introduce the following notation
\be\la{vorm}\l\{\ba{l}a_0(x):=a_0(n,\a,x),\;\ndots\\ \\
a_i(x):=a_i(\a,x),\;i=1,2,3,\ldots.\ea\r.\ee

In order to find the general form of the coefficients $\l\{a_i(x)\r\}_{i=0}^{\infty}$
we will prove the following theorems.

\st{The polynomials $\set{\SGP}$ satisfy the infinite order differential
equation given by
\be\la{dv1}\sum_{i=0}^{\infty}b_i(x)y^{(i)}(x)=0\ee
and
\be\la{b}\l\{\ba{l}\ds b_0(x):=b_0(n,\a,x)=\frac{1}{2}\l[1-(-1)^n\r],\;\ndots\\
\\
\ds b_i(x):=b_i(\a,x)=\frac{2^{i-1}}{i!}(-x)^i,\;i=1,2,3,\ldots.\ea\r.\ee}

\noindent
Note that the coefficients $\l\{b_i(x)\r\}_{i=0}^{\infty}$ do not depend
on $\a$. The proof of theorem 1 can be found in the next section.

\st{The polynomials $\set{\SGP}$ satisfy the differential equation given by
\be\la{dv2}M\sum_{i=0}^{\infty}c_i(x)y^{(i)}(x)+(1-x^2)y''(x)-2(\a+1)xy'(x)+n(n+2\a+1)y(x)=0\ee
where the coefficients $\l\{c_i(x)\r\}_{i=0}^{\infty}$ are defined by
\bea\la{cnul}c_0(x):=c_0(n,\a,x)=4(2\a+3){n+2\a+2 \ch n-2},\;\ndots\eea
and
\be\la{cstar}c_i(x)=(2\a+3)(1-x^2)c_i^*(x),\;i=1,2,3,\ldots,\ee
where
\be\la{c}\l\{\ba{l}\ds c_1^*(x):=c_1^*(\a,x)=0\\ \\
\ds c_i^*(x):=c_i^*(\a,x)=\frac{2^i}{i!}\sum_{k=0}^{i-2}
{\a+1 \ch i-k-2}{i-2\a-5 \ch k}\l(\frac{1-x}{2}\r)^k,\\
\hfill i=2,3,4,\ldots.\ea\r.\ee}

\noindent
The proof of theorem 2 will be given in section 5.

Now we will show that the general solution $\l\{a_i(x)\r\}_{i=0}^{\infty}$ of
(\ref{rel1}) and (\ref{rel2}), where $\l\{a_i(x)\r\}_{i=0}^{\infty}$ are
continuous functions on the real line and where $\l\{a_i(x)\r\}_{i=1}^{\infty}$
are independent of $n$, is
\be\la{a}\l\{\ba{l}a_0(n,\a,x)=a_0(1,\a,x)b_0(n,\a,x)+c_0(n,\a,x),\;\ndots\\ \\
a_i(\a,x)=a_0(1,\a,x)b_i(\a,x)+c_i(\a,x),\;i=1,2,3,\ldots,\ea\r.\ee
where $a_0(1,\a,x)$ is an arbitrary continuous function on the real line and
where the coefficients $\l\{b_i(x)\r\}_{i=0}^{\infty}$ and
$\l\{c_i(x)\r\}_{i=0}^{\infty}$ are given by (\ref{b}), (\ref{cnul}),
(\ref{cstar}) and (\ref{c}).

The proof is based on the following lemma.

\par\vspace{5mm}
{\bf Lemma.} {\sl Let $\l\{a_i(x)\r\}_{i=0}^{\infty}$ be a set of continuous
functions of the form (\ref{vorm}) which satisfies the homogeneous system
\be\la{hom}\l\{\ba{l}\ds\sum_{i=0}^{\infty}a_i(x)D^i\SP=0\\
\\
\ds\sum_{i=0}^{\infty}ia_i(x)D^i\SP+x\sum_{i=0}^{\infty}a_i(x)D^{i+1}\SP=0
\ea\r.\ee
with $a_0(1,\a,x)=0$ for all real $x$. Then we have
$$a_i(x)=0,\;i=0,1,2,\ldots$$
for all real $x$.}
\par\vspace{5mm}

In order to prove this lemma we substitute small values of $n$ in
the homogeneous system (\ref{hom}). Then we obtain for $n=0$ and $n=1$~:
$$a_0(0,\a,x)=0\;\;\mbox{ and }\;\;(\a+1)\l[xa_0(1,\a,x)+a_1(\a,x)\r]=0$$
for all real $x$. Since $a_0(1,\a,x)=0$ and $\a>-1$, we conclude that
$a_1(\a,x)=0$ for all real $x$. If we substitute $n=2$ in (\ref{hom}) we
obtain
$$\l\{\ba{l}\ds a_0(2,\a,x)P_2^{(\a,\a)}(x)+
a_2(\a,x)\frac{d^2}{dx^2}P_2^{(\a,\a)}(x)=0\\ \\
\ds 2a_2(\a,x)\frac{d^2}{dx^2}P_2^{(\a,\a)}(x)+
xa_0(2,\a,x)\frac{d}{dx}P_2^{(\a,\a)}(x)=0.\ea\r.$$
Hence
$$\l\{\ba{l}\ds \l[2P_2^{(\a,\a)}(x)-x\frac{d}{dx}P_2^{(\a,\a)}(x)\r]
a_0(2,\a,x)=0\\ \\
\ds a_2(\a,x)\frac{d^2}{dx^2}P_2^{(\a,\a)}(x)=-a_0(2,\a,x)P_2^{(\a,\a)}(x).
\ea\r.$$
Since
$$2P_2^{(\a,\a)}(x)-x\frac{d}{dx}P_2^{(\a,\a)}(x)=-\frac{1}{2}(\a+2)\ne 0$$
we conclude that $a_0(2,\a,x)=0$ for all real $x$ and therefore that
$a_2(\a,x)=0$ for all real $x$, since
$$\frac{d^2}{dx^2}P_2^{(\a,\a)}(x)=\frac{1}{2}(\a+2)(2\a+3)\ne 0.$$
In the same way we obtain for $n=3$~:
$$\l\{\ba{l}\ds \l[3P_3^{(\a,\a)}(x)-x\frac{d}{dx}P_3^{(\a,\a)}(x)\r]
a_0(3,\a,x)=0\\ \\
\ds a_3(\a,x)\frac{d^3}{dx^3}P_3^{(\a,\a)}(x)=-a_0(3,\a,x)P_3^{(\a,\a)}(x).
\ea\r.$$
Now we find $a_0(3,\a,x)=0$ for all real $x$ except for $x=0$ being the only
zero of
$$3P_3^{(\a,\a)}(x)-x\frac{d}{dx}P_3^{(\a,\a)}(x)=-\frac{1}{2}(\a+2)(\a+3)x.$$
The continuity of $a_0(3,\a,x)$ implies that $a_0(3,\a,x)=0$ for all real $x$.
Then we also have $a_3(\a,x)=0$ for all real $x$, since
$$\frac{d^3}{dx^3}P_3^{(\a,\a)}(x)=\frac{1}{2}(\a+2)(\a+3)(2\a+5)\ne 0.$$
If we proceed in this way we also find for each $n\ge 4$ that $a_0(n,\a,x)=0$
for all real $x$ except for the possible zeros of $n\SP-x\frac{d}{dx}\SP$.
The continuity of $a_0(n,\a,x)$ then implies that $a_0(n,\a,x)=0$ for all
real $x$ and finally we have $a_n(\a,x)=0$ for all real $x$, since
$$\frac{d^n}{dx^n}\SP={2n+2\a \ch n}\frac{n!}{2^n}\ne 0.$$
This completes the proof of the lemma.

For the moment we introduce the term {\em continuous sequence}. A sequence
$\l\{a_i(x)\r\}_{i=0}^{\infty}$ is called a continuous sequence if
$\l\{a_i(x)\r\}_{i=0}^{\infty}$ are continuous functions on the real line
of the form (\ref{vorm}).

In section 4 (proof of theorem 1) we show that the continuous sequence
$\l\{b_i(x)\r\}_{i=0}^{\infty}$ defined by (\ref{b}) is a solution of the
homogeneous system (\ref{hom}).
This implies that the general solution $\l\{a_i(x)\r\}_{i=0}^{\infty}$ of
(\ref{hom}), being a continuous sequence, is
$$a_i(x)=a_0(1,\a,x)b_i(x),\;i=0,1,2,\ldots.$$
This can be shown by using the lemma as follows. We define
$$a_i^*(x):=a_i(x)-a_0(1,\a,x)b_i(x),\;i=0,1,2,\ldots.$$
Since $a_0(1,\a,x)$ is continuous, $\l\{a_i^*(x)\r\}_{i=0}^{\infty}$ is a
continuous sequence which is a solution of (\ref{hom}) too where
$$a_0^*(1,\a,x)=a_0(1,\a,x)-a_0(1,\a,x)b_0(1,\a,x)=0$$
for all real $x$. Now the lemma gives us the desired result.

In section 5 (proof of theorem 2) we show that the continuous sequence
$\l\{c_i(x)\r\}_{i=0}^{\infty}$ given by (\ref{cnul}), (\ref{cstar}) and
(\ref{c}) is a solution of (\ref{rel1}) and (\ref{rel2}).
In order to prove that the general solution
$\l\{a_i(x)\r\}_{i=0}^{\infty}$ of (\ref{rel1}) and (\ref{rel2}), being a
continuous sequence, is of the form (\ref{a}), we simply note that
$\l\{a_i(x)-c_i(x)\r\}_{i=0}^{\infty}$ is a continuous sequence and a
solution of the homogeneous system (\ref{hom}). Hence
$$a_i(x)-c_i(x)=a_0(1,\a,x)b_i(x),\;i=0,1,2,\ldots.$$

This shows that we have found all differential equations of the form (\ref{dv}),
where the coefficients $\l\{a_i(x)\r\}_{i=0}^{\infty}$ are continuous functions
on the real line and where $\l\{a_i(x)\r\}_{i=1}^{\infty}$ are independent of
$n$.

In \cite{Gran} we conjectured that the polynomials $\set{\SGP}$ satisfy a
differential equation of the form (\ref{dv2}) where the coefficients
$\l\{c_i(x)\r\}_{i=0}^{\infty}$ are given by (\ref{cnul}) and (\ref{cstar})
where
\be\la{cgran}\l\{\ba{l}\ds c_{2i}^*(\a,x)=\frac{4(-1)^{i+1}}{(2i)!}{\a+1 \ch i-1}
\hyp{2}{1}{-i+1,\a+\frac{5}{2}-i}{\frac{1}{2}}{x^2},\;i=1,2,3,\ldots\\ \\
\ds c_1^*(\a,x)=0\\ \\
\ds c_{2i+1}^*(\a,x)=\frac{8(-1)^{i+1}}{(2i+1)!}{\a \ch i-1}(\a+1)x\;
\hyp{2}{1}{-i+1,\a+\frac{5}{2}-i}{\frac{3}{2}}{x^2},\\
\hfill i=1,2,3,\ldots.\ea\r.\ee
By using the formulas (\ref{even}) and (\ref{odd}) we can prove that
(\ref{cgran}) is equivalent to
\be\la{equiv}\l\{\ba{l}\ds c_1^*(\a,x)=0\\ \\
\ds c_i^*(\a,x)=\frac{2^i}{i!}P_{i-2}^{(\a-i+3,\a-i+3)}(x),
\;i=2,3,4,\ldots.\ea\r.\ee
Here we remark that the parameter $\a-i+3$ in the Jacobi polynomial might be
smaller than $-1$. However, the ultraspherical polynomial $\SP$ is also a
polynomial in $\a$. Instead of the definition (\ref{defJ}) we may define
$$\SP=\frac{1}{n!}\sum_{k=0}^n{n \ch k}(n+2\a+1)_k(\a+k+1)_{n-k}\l(\frac{x-1}{2}\r)^k,
\;\ndots.$$
By using this definition of the ultraspherical polynomial we see that (\ref{equiv})
is equivalent to (\ref{c}).

Since, by using (\ref{a}), (\ref{b}), (\ref{cstar}) and (\ref{cgran}) we have
$$a_{2i}(0):=a_{2i}(\a,0)=
(2\a+3)\frac{4(-1)^{i+1}}{(2i)!}{\a+1 \ch i-1}\ne 0,\;i=1,2,3,\ldots$$
if $\a$ is not a nonnegative integer we conclude that if $M>0$ the order of
the differential equation given by (\ref{dv}) can only be finite for
nonnegative integer values of $\a$.
From the form (\ref{cgran}) we easily see that for nonnegative integer
values of $\a$ we have
$$\l\{\ba{l}\mbox{degree}[c_i(x)]=i,\;i=2,3,4,\ldots,2\a+4\\ \\
c_i(x)=0,\;i>2\a+4.\ea\r.$$
Hence, for $M>0$ the order of the differential equation given by (\ref{dv})
can only be finite if we choose $a_0(1,\a,x)=0$ in view of (\ref{a}) and
if $\a$ is a nonnegative integer. In that case the order of the differential
equation given by (\ref{dv2}), (\ref{cnul}), (\ref{cstar}) and (\ref{c})
equals $2\a+4$.

\section{Proof of theorem 1.}

To prove theorem 1 we have to show that
$$\l\{\ba{l}\ds \sum_{i=0}^{\infty}b_i(x)D^i\SP=0\\ \\
\ds \sum_{i=0}^{\infty}ib_i(x)D^i\SP+x\sum_{i=0}^{\infty}b_i(x)D^{i+1}\SP=0.\ea\r.$$
To do this we first note that we have in view of (\ref{b})
$$\sum_{i=0}^{\infty}b_i(x)D^i\SP=\frac{1}{2}\l[1-(-1)^n\r]\SP+
\sum_{i=1}^{\infty}b_i(x)D^i\SP$$
and
$$\sum_{i=0}^{\infty}b_i(x)D^{i+1}\SP=\frac{1}{2}\l[1-(-1)^n\r]\frac{d}{dx}\SP
+\sum_{i=1}^{\infty}b_i(x)D^{i+1}\SP.$$
Now we use (\ref{afgJ}), (\ref{symJ}) and (\ref{b}) to obtain
\bea\sum_{i=1}^{\infty}b_i(x)D^i\SP
&=&\frac{1}{2}{n+\a \ch n}\sum_{i=1}^{\infty}\frac{x^i}{i!}\sum_{k=i}^{\infty}
\frac{(-n)_k(n+2\a+1)_k}{(\a+1)_k(k-i)!}\l(\frac{1-x}{2}\r)^{k-i}\nn
&=&\frac{1}{2}{n+\a \ch n}\sum_{k=1}^{\infty}\frac{(-n)_k(n+2\a+1)_k}{(\a+1)_kk!}
\sum_{i=1}^k{k \ch i}x^i\l(\frac{1-x}{2}\r)^{k-i}\nn
&=&\frac{1}{2}{n+\a \ch n}\sum_{k=1}^{\infty}\frac{(-n)_k(n+2\a+1)_k}{(\a+1)_kk!}
\l[\l(\frac{1+x}{2}\r)^k-\l(\frac{1-x}{2}\r)^k\r]\nn
&=&\frac{1}{2}\l[P_n^{(\a,\a)}(-x)-\SP\r]=-\frac{1}{2}\l[1-(-1)^n\r]\SP,\n\eea
\bea\sum_{i=0}^{\infty}ib_i(x)D^i\SP
&=&\frac{1}{2}{n+\a \ch n}\sum_{i=1}^{\infty}\frac{x^i}{(i-1)!}\sum_{k=i}^{\infty}
\frac{(-n)_k(n+2\a+1)_k}{(\a+1)_k(k-i)!}\l(\frac{1-x}{2}\r)^{k-i}\nn
&=&\frac{1}{2}{n+\a \ch n}\sum_{k=1}^{\infty}\frac{(-n)_k(n+2\a+1)_k}{(\a+1)_k(k-1)!}
\sum_{i=1}^k{k-1 \ch i-1}x^i\l(\frac{1-x}{2}\r)^{k-i}\nn
&=&\frac{1}{2}{n+\a \ch n}x\sum_{k=1}^{\infty}\frac{(-n)_k(n+2\a+1)_k}{(\a+1)_k(k-1)!}
\l(\frac{1+x}{2}\r)^{k-1}\nn
&=&(-1)^nx\frac{d}{dx}\SP,\n\eea
and
\bea & &\sum_{i=1}^{\infty}b_i(x)D^{i+1}\SP\nn
&=&-\frac{1}{4}{n+\a \ch n}\sum_{i=1}^{\infty}\frac{x^i}{i!}\sum_{k=i+1}^{\infty}
\frac{(-n)_k(n+2\a+1)_k}{(\a+1)_k(k-i-1)!}\l(\frac{1-x}{2}\r)^{k-i-1}\nn
&=&-\frac{1}{4}{n+\a \ch n}\sum_{k=2}^{\infty}\frac{(-n)_k(n+2\a+1)_k}
{(\a+1)_k(k-1)!}\sum_{i=1}^{k-1}{k-1 \ch i}x^i\l(\frac{1-x}{2}\r)^{k-i-1}\nn
&=&-\frac{1}{4}{n+\a \ch n}\sum_{k=2}^{\infty}\frac{(-n)_k(n+2\a+1)_k}
{(\a+1)_k(k-1)!}\l[\l(\frac{1+x}{2}\r)^{k-1}-\l(\frac{1-x}{2}\r)^{k-1}\r]\nn
&=&-\frac{1}{2}\l[1+(-1)^n\r]\frac{d}{dx}\SP.\n\eea
This proves theorem 1.

\section{Proof of theorem 2.}

In this section we will give a proof of the main theorem 2.

In view of (\ref{rel1}) and (\ref{rel2}) we have to show that
\be\la{first}\sum_{i=0}^{\infty}c_i(x)D^i\SP=\frac{4}{(2\a+1)}{n+2\a \ch n}
\frac{d^2}{dx^2}\SP\ee
and
\be\la{second}\sum_{i=0}^{\infty}ic_i(x)D^i\SP+x\sum_{i=0}^{\infty}c_i(x)D^{i+1}\SP
=4{n+2\a+1 \ch n-1}\frac{d^2}{dx^2}\SP.\ee

Now we write by using (\ref{cnul}), (\ref{cstar}) and (\ref{c})
$$\sum_{i=0}^{\infty}c_i(x)D^i\SP=4(2\a+3){n+2\a+2 \ch n-2}\SP+
\sum_{i=2}^{\infty}c_i(x)D^i\SP,$$
$$\sum_{i=0}^{\infty}ic_i(x)D^i\SP=\sum_{i=2}^{\infty}ic_i(x)D^i\SP$$
and
$$\sum_{i=0}^{\infty}c_i(x)D^{i+1}\SP=4(2\a+3){n+2\a+2 \ch n-2}\frac{d}{dx}\SP
+\sum_{i=2}^{\infty}c_i(x)D^{i+1}\SP.$$

Now we obtain by using (\ref{c})
\bea & &\sum_{i=2}^{\infty}c_i^*(x)D^i\SP\nn
&=&\sum_{i=2}^{\infty}\frac{2^i}{i!}\sum_{k=0}^{i-2}{\a+1 \ch i-k-2}
{i-2\a-5 \ch k}\l(\frac{1-x}{2}\r)^kD^i\SP\nn
&=&\sum_{k=0}^{\infty}\sum_{i=k+2}^{\infty}\frac{2^i}{i!}{\a+1 \ch i-k-2}
{i-2\a-5 \ch k}\l(\frac{1-x}{2}\r)^kD^i\SP\nn
&=&\sum_{k=0}^{\infty}\sum_{i=0}^{\infty}\frac{2^{i+k+2}}{(i+k+2)!}
{\a+1 \ch i}{i+k-2\a-3 \ch k}\l(\frac{1-x}{2}\r)^kD^{i+k+2}\SP\nn
&=&\sum_{i=0}^{\infty}{\a+1 \ch i}\frac{2^{i+2}}{(i+2)!}\sum_{k=0}^{\infty}
\frac{(-2\a+i-2)_k}{k!(i+3)_k}(1-x)^kD^{i+k+2}\SP.\n\eea
By using (\ref{afgJ}) we find
\bea & &D^{i+k+2}\SP\nn
&=&{n+\a \ch n}\l(-\frac{1}{2}\r)^{i+k+2}\sum_{m=i+k+2}^{\infty}
\frac{(-n)_m(n+2\a+1)_m}{(m-i-k-2)!(\a+1)_m}\l(\frac{1-x}{2}\r)^{m-i-k-2}\nn
&=&{n+\a \ch n}\l(-\frac{1}{2}\r)^{i+k+2}\sum_{m=k}^{\infty}
\frac{(-n)_{m+i+2}(n+2\a+1)_{m+i+2}}{(m-k)!(\a+1)_{m+i+2}}
\l(\frac{1-x}{2}\r)^{m-k}.\n\eea
Now we use Vandermonde's summation formula (\ref{Van}) to see that this leads to
\bea & &\sum_{k=0}^{\infty}\frac{(-2\a+i-2)_k}{k!(i+3)_k}(1-x)^kD^{i+k+2}\SP\nn
&=&{n+\a \ch n}\sum_{k=0}^{\infty}\sum_{m=k}^{\infty}\frac{(-2\a+i-2)_k}{k!(i+3)_k}
\frac{(-n)_{m+i+2}(n+2\a+1)_{m+i+2}}{(m-k)!(\a+1)_{m+i+2}}(-1)^{i+k}\frac{(1-x)^m}{2^{i+m+2}}\nn
&=&{n+\a \ch n}\sum_{m=0}^{\infty}\sum_{k=0}^m\frac{(-2\a+i-2)_k}{k!(i+3)_k}
\frac{(-n)_{m+i+2}(n+2\a+1)_{m+i+2}}{2^{m+i+2}(\a+1)_{m+i+2}}\frac{(-m)_k}{m!}(-1)^i(1-x)^m\nn
&=&{n+\a \ch n}(-1)^i\sum_{m=0}^{\infty}\frac{(-n)_{m+i+2}(n+2\a+1)_{m+i+2}}
{2^{i+m+2}(\a+1)_{m+i+2}}\frac{(1-x)^m}{m!}\times{}\nn
& &\hspace{8cm}{}\times\hyp{2}{1}{-m,-2\a+i-2}{i+3}{1}\nn
&=&{n+\a \ch n}(-1)^i\sum_{m=0}^{\infty}\frac{(-n)_{m+i+2}(n+2\a+1)_{m+i+2}}
{2^{m+i+2}(\a+1)_{m+i+2}}\frac{(1-x)^m}{m!}\frac{(2\a+5)_m}{(i+3)_m}\nn
&=&{n+\a \ch n}(-1)^i\frac{(i+2)!}{2^{i+2}}\sum_{m=0}^{\infty}\frac{(-n)_{m+i+2}
(n+2\a+1)_{m+i+2}}{(m+i+2)!(\a+1)_{m+i+2}}\frac{(2\a+5)_m}{m!}\l(\frac{1-x}{2}\r)^m.\n\eea
Hence
\bea & &\sum_{i=2}^{\infty}c_i^*(x)D^i\SP\nn
&=&\sum_{i=0}^{\infty}{\a+1 \ch i}\frac{2^{i+2}}{(i+2)!}\sum_{k=0}^{\infty}
\frac{(-2\a+i-2)_k}{k!(i+3)_k}(1-x)^kD^{i+k+2}\SP\nn
&=&{n+\a \ch n}\sum_{i=0}^{\infty}(-1)^i{\a+1 \ch i}\sum_{m=0}^{\infty}\frac{(-n)_{m+i+2}
(n+2\a+1)_{m+i+2}}{(m+i+2)!(\a+1)_{m+i+2}}\frac{(2\a+5)_m}{m!}\l(\frac{1-x}{2}\r)^m\nn
&=&{n+\a \ch n}\sum_{m=0}^{\infty}\frac{(-n)_{m+2}(n+2\a+1)_{m+2}}{(m+2)!(\a+1)_{m+2}}
\frac{(2\a+5)_m}{m!}\l(\frac{1-x}{2}\r)^m\times{}\nn
& &\hspace{5cm}{}\times\hyp{3}{2}{-n+m+2,-\a-1,n+2\a+m+3}{m+3,\a+m+3}{1}.\n\eea
Now we use the well-known fact that
$$\hyp{3}{2}{a,b,c}{d,e}{z}=\frac{c}{c-a}\hyp{3}{2}{a,b,c+1}{d,e}{z}
-\frac{a}{c-a}\hyp{3}{2}{a+1,b,c}{d,e}{z},\;c-a\ne 0$$
and the Saalsch\"{u}tz summation formula (\ref{Saal}) to find for $n>m+2$
\bea & &\hyp{3}{2}{-n+m+2,-\a-1,n+2\a+m+3}{m+3,\a+m+3}{1}\nn
&=&\frac{(n+2\a+m+3)}{(2n+2\a+1)}\hyp{3}{2}{-n+m+2,-\a-1,n+2\a+m+4}{m+3,\a+m+3}{1}+{}\nn
& &\hspace{2cm}{}+\frac{(n-m-2)}{(2n+2\a+1)}\hyp{3}{2}{-n+m+3,-\a-1,n+2\a+m+3}{m+3,\a+m+3}{1}\nn
&=&\frac{(n+2\a+m+3)}{(2n+2\a+1)}\frac{(\a+m+4)_{n-m-2}(-n-2\a-1)_{n-m-2}}{(m+3)_{n-m-2}(-n-\a)_{n-m-2}}+{}\nn
& &\hspace{2cm}{}+\frac{(n-m-2)}{(2n+2\a+1)}\frac{(\a+m+4)_{n-m-3}(-n-2\a)_{n-m-3}}{(m+3)_{n-m-3}(-n-\a+1)_{n-m-3}}\nn
&=&\frac{\G(n+2\a+1)}{n!}\frac{1}{(\a+m+3)}\frac{(m+2)!}{\G(m+2\a+4)}
\l[n(n+2\a+1)+(\a+1)(2\a+m+3)\r].\n\eea
Note that the same result also holds for $n=m+2$.
For $n<m+2$ we have $(-n)_{m+2}=0$ and the $_3F_2(1)$ exists since
$$(-n+m+2)+(-\a-1)+(n+2\a+m+3)<(m+3)+(\a+m+3).$$
This implies that
\bea & &\sum_{i=2}^{\infty}c_i^*(x)D^i\SP\nn
&=&{n+\a \ch n}\frac{\G(n+2\a+1)}{n!}\sum_{m=0}^{\infty}\frac{(-n)_{m+2}
(n+2\a+1)_{m+2}}{m!(\a+1)_{m+3}}\frac{(2\a+5)_m}{\G(m+2\a+4)}\times{}\nn
& &\hspace{4cm}{}\times\l[n(n+2\a+1)+(\a+1)(2\a+m+3)\r]\l(\frac{1-x}{2}\r)^m\nn
&=&\frac{1}{(2\a+1)_4}{n+\a \ch n}{n+2\a \ch n}\sum_{m=0}^{\infty}
\frac{(-n)_{m+2}(n+2\a+1)_{m+2}}{m!(\a+1)_{m+3}}(m+2\a+4)\times{}\nn
& &\hspace{4cm}{}\times\l[n(n+2\a+1)+(\a+1)(m+2\a+3)\r]\l(\frac{1-x}{2}\r)^m.\n\eea
Since we have
$$n(n+2\a+1)+(\a+1)(m+2\a+3)=(\a+m+3)(m+2\a+3)+(n-m-2)(n+2\a+m+3)$$
and
$$(m+2\a+3)(m+2\a+4)=m(m-1)+2(2\a+4)m+(2\a+3)(2\a+4)$$
we obtain
\bea & &\sum_{i=2}^{\infty}c_i^*(x)D^i\SP\nn
&=&\frac{1}{(2\a+1)_4}{n+\a \ch n}{n+2\a \ch n}\times{}\nn
& &\hspace{1cm}{}\times\l[\sum_{m=0}^{\infty}\frac{(-n)_{m+4}(n+2\a+1)_{m+4}}
{m!(\a+1)_{m+4}}\l(\frac{1-x}{2}\r)^{m+2}+{}\r.\nn
& &\hspace{2cm}{}+2(2\a+4)\sum_{m=0}^{\infty}
\frac{(-n)_{m+3}(n+2\a+1)_{m+3}}{m!(\a+1)_{m+3}}\l(\frac{1-x}{2}\r)^{m+1}+{}\nn
& &\hspace{2cm}{}+(2\a+3)(2\a+4)\sum_{m=0}^{\infty}
\frac{(-n)_{m+2}(n+2\a+1)_{m+2}}{m!(\a+1)_{m+2}}\l(\frac{1-x}{2}\r)^m+{}\nn
& &\hspace{2cm}{}-\sum_{m=0}^{\infty}\frac{(-n)_{m+4}(n+2\a+1)_{m+4}}
{m!(\a+1)_{m+4}}\l(\frac{1-x}{2}\r)^{m+1}+{}\nn
& &\hspace{2cm}\l.{}-(2\a+4)\sum_{m=0}^{\infty}
\frac{(-n)_{m+3}(n+2\a+1)_{m+3}}{m!(\a+1)_{m+3}}\l(\frac{1-x}{2}\r)^m\r]\nn
&=&\frac{4}{(2\a+1)_4}{n+2\a \ch n}\times{}\nn
& &\hspace{1cm}{}\times\l[4\l(\frac{1-x}{2}\r)^2\frac{d^4}{dx^4}\SP-4(2\a+4)
\l(\frac{1-x}{2}\r)\frac{d^3}{dx^3}\SP+{}\r.\nn
& &\hspace{2cm}{}+(2\a+3)(2\a+4)\frac{d^2}{dx^2}\SP+{}\nn
& &\hspace{2cm}\l.{}-4\l(\frac{1-x}{2}\r)\frac{d^4}{dx^4}\SP+2(2\a+4)\frac{d^3}{dx^3}\SP\r]\nn
&=&\frac{4}{(2\a+1)_4}{n+2\a \ch n}\times{}\nn
& &\hspace{1cm}{}\times\l[-(1-x^2)\frac{d^4}{dx^4}\SP+2(2\a+4)x\frac{d^3}{dx^3}\SP+{}\r.\nn
& &\hspace{6cm}\l.{}+(2\a+3)(2\a+4)\frac{d^2}{dx^2}\SP\r].\n\eea
Now we use (\ref{dvJ}) for $i=2$ to find
\bea & &\sum_{i=2}^{\infty}c_i^*(x)D^i\SP\nn
&=&\frac{4}{(2\a+1)_4}{n+2\a \ch n}\times{}\nn
& &{}\times\l[2(\a+1)x\frac{d^3}{dx^3}\SP
+\l[(n-2)(n+2\a+3)+(2\a+3)(2\a+4)\r]\frac{d^2}{dx^2}\SP\r].\n\eea
Hence, by using (\ref{cstar}) and (\ref{dvJ}) for $i=1$ and $i=0$ we find
\bea & &\sum_{i=2}^{\infty}c_i(x)D^i\SP\nn
&=&\frac{1}{(\a+1)(\a+2)(2\a+1)}{n+2\a \ch n}(1-x^2)\times{}\nn
& &{}\times\l[2(\a+1)x\frac{d^3}{dx^3}\SP
+\l[(n-2)(n+2\a+3)+(2\a+3)(2\a+4)\r]\frac{d^2}{dx^2}\SP\r]\nn
&=&\frac{1}{(\a+1)(\a+2)(2\a+1)}{n+2\a \ch n}\times{}\nn
& &{}\times\l[4(\a+1)(\a+2)x^2\frac{d^2}{dx^2}\SP
-(n-1)(n+2\a+2)(1-x^2)\frac{d^2}{dx^2}\SP+{}\r.\nn
& &\hspace{2cm}{}-n(n-1)(n+2\a+1)(n+2\a+2)\SP+{}\nn
& &\hspace{2cm}\l.{}+\l[(n-2)(n+2\a+3)+(2\a+3)(2\a+4)\r](1-x^2)\frac{d^2}{dx^2}\SP\r]\nn
&=&\frac{4}{(2\a+1)}{n+2\a \ch n}\frac{d^2}{dx^2}\SP+{}\nn
& &\hspace{2cm}{}-\frac{n(n-1)(n+2\a+1)(n+2\a+2)}{(\a+1)(\a+2)(2\a+1)}{n+2\a \ch n}\SP\nn
&=&\frac{4}{(2\a+1)}{n+2\a \ch n}\frac{d^2}{dx^2}\SP-
4(2\a+3){n+2\a+2 \ch n-2}\SP.\n\eea
Finally, this implies that
$$\sum_{i=0}^{\infty}c_i(x)D^i\SP=\frac{4}{(2\a+1)}{n+2\a \ch n}\frac{d^2}{dx^2}\SP,$$
which proves (\ref{first}).

The proof of (\ref{second}) is much easier. We start from
\bea & &\sum_{i=2}^{\infty}ic_i^*(x)D^i\SP\nn
&=&\sum_{i=0}^{\infty}{\a+1 \ch i}\frac{2^{i+2}}{(i+1)!}\sum_{k=0}^{\infty}
\frac{(-2\a+i-2)_k}{k!(i+2)_k}(1-x)^kD^{i+k+2}\SP.\n\eea
We use Vandermonde's summation formula (\ref{Van}) again to obtain
\bea & &\sum_{k=0}^{\infty}\frac{(-2\a+i-2)_k}{k!(i+2)_k}(1-x)^kD^{i+k+2}\SP\nn
&=&{n+\a \ch n}\sum_{m=0}^{\infty}\sum_{k=0}^m\frac{(-2\a+i-2)_k}{k!(i+2)_k}
\frac{(-n)_{m+i+2}(n+2\a+1)_{m+i+2}}{2^{m+i+2}(\a+1)_{m+i+2}}\frac{(-m)_k}{m!}
(-1)^i(1-x)^m\nn
&=&{n+\a \ch n}(-1)^i\sum_{m=0}^{\infty}\frac{(-n)_{m+i+2}(n+2\a+1)_{m+i+2}}
{2^{m+i+2}(\a+1)_{m+i+2}}\frac{(1-x)^m}{m!}\times{}\nn
& &\hspace{8cm}{}\times\hyp{2}{1}{-m,-2\a+i-2}{i+2}{1}\nn
&=&{n+\a \ch n}(-1)^i\sum_{m=0}^{\infty}\frac{(-n)_{m+i+2}(n+2\a+1)_{m+i+2}}
{2^{m+i+2}(\a+1)_{m+i+2}}\frac{(1-x)^m}{m!}\frac{(2\a+4)_m}{(i+2)_m}\nn
&=&{n+\a \ch n}(-1)^i\frac{(i+1)!}{2^{i+2}}\sum_{m=0}^{\infty}\frac{(-n)_{m+i+2}
(n+2\a+1)_{m+i+2}}{(m+i+1)!(\a+1)_{m+i+2}}\frac{(2\a+4)_m}{m!}\l(\frac{1-x}{2}\r)^m.\n\eea
This leads to
\bea & &\sum_{i=2}^{\infty}ic_i^*(x)D^i\SP\nn
&=&{n+\a \ch n}\sum_{i=0}^{\infty}(-1)^i{\a+1 \ch i}\sum_{m=0}^{\infty}
\frac{(-n)_{m+i+2}(n+2\a+1)_{m+i+2}}{(m+i+1)!(\a+1)_{m+i+2}}\frac{(2\a+4)_m}{m!}
\l(\frac{1-x}{2}\r)^m\nn
&=&{n+\a \ch n}\sum_{m=0}^{\infty}\frac{(-n)_{m+2}(n+2\a+1)_{m+2}}{(m+1)!(\a+1)_{m+2}}
\frac{(2\a+4)_m}{m!}\l(\frac{1-x}{2}\r)^m\times{}\nn
& &\hspace{5cm}{}\times\hyp{3}{2}{-n+m+2,-\a-1,n+2\a+m+3}{m+2,\a+m+3}{1}.\n\eea
Now we use the Saalsch\"{u}tz summation formula (\ref{Saal}) again to find
for $n\ge m+2$
\bea & &\hyp{3}{2}{-n+m+2,-\a-1,n+2\a+m+3}{m+2,\a+m+3}{1}\nn
&=&\frac{(\a+m+3)_{n-m-2}(-n-2\a-1)_{n-m-2}}{(m+2)_{n-m-2}(-n-\a)_{n-m-2}}\nn
&=&\frac{\G(n+\a+1)}{\G(\a+m+3)}\frac{(m+1)!}{\G(n)}\frac{\G(n+2\a+2)}{\G(m+2\a+4)}
\frac{\G(\a+m+3)}{\G(n+\a+1)}=\frac{(m+1)!}{\G(n)}\frac{\G(n+2\a+2)}{\G(m+2\a+4)}.\n\eea
Since $(-n)_{m+2}=0$ and the $_3F_2(1)$ exists for $n<m+2$,
this implies, by using (\ref{afgJ})
\bea & &\sum_{i=2}^{\infty}ic_i^*(x)D^i\SP\nn
&=&{n+\a \ch n}\sum_{m=0}^{\infty}\frac{(-n)_{m+2}(n+2\a+1)_{m+2}}{(m+1)!(\a+1)_{m+2}}
\frac{(2\a+4)_m}{m!}\frac{(m+1)!}{\G(n)}\frac{\G(n+2\a+2)}{\G(m+2\a+4)}\l(\frac{1-x}{2}\r)^m\nn
&=&{n+\a \ch n}\frac{\G(n+2\a+2)}{\G(n)\G(2\a+4)}\sum_{m=0}^{\infty}
\frac{(-n)_{m+2}(n+2\a+1)_{m+2}}{m!(\a+1)_{m+2}}\l(\frac{1-x}{2}\r)^m\nn
&=&\frac{4}{(2\a+3)}{n+2\a+1 \ch n-1}\frac{d^2}{dx^2}\SP,\n\eea
and therefore we have
$$\sum_{i=2}^{\infty}ic_i(x)D^i\SP=4{n+2\a+1 \ch n-1}(1-x^2)\frac{d^2}{dx^2}\SP.$$

Now we look at the second sum on the left-hand side of (\ref{second}).
Now we have
\bea & &\sum_{i=2}^{\infty}c_i^*(x)D^{i+1}\SP\nn
&=&\sum_{i=0}^{\infty}{\a+1 \ch i}\frac{2^{i+2}}{(i+2)!}\sum_{k=0}^{\infty}
\frac{(-2\a+i-2)_k}{k!(i+3)_k}(1-x)^kD^{i+k+3}\SP.\n\eea
As before we find by using the Vandermonde summation formula (\ref{Van})
\bea & &\sum_{k=0}^{\infty}\frac{(-2\a+i-2)_k}{k!(i+3)_k}(1-x)^kD^{i+k+3}\SP\nn
&=&{n+\a \ch n}\sum_{m=0}^{\infty}\sum_{k=0}^m\frac{(-2\a+i-2)_k}{k!(i+3)_k}
\frac{(-n)_{m+i+3}(n+2\a+1)_{m+i+3}}{2^{m+i+3}(\a+1)_{m+i+3}}\frac{(-m)_k}{m!}
(-1)^{i+1}(1-x)^m\nn
&=&{n+\a \ch n}(-1)^{i+1}\sum_{m=0}^{\infty}\frac{(-n)_{m+i+3}(n+2\a+1)_{m+i+3}}
{2^{m+i+3}(\a+1)_{m+i+3}}\frac{(1-x)^m}{m!}\times{}\nn
& &\hspace{8cm}{}\times\hyp{2}{1}{-m,-2\a+i-2}{i+3}{1}\nn
&=&{n+\a \ch n}(-1)^{i+1}\sum_{m=0}^{\infty}\frac{(-n)_{m+i+3}(n+2\a+1)_{m+i+3}}
{2^{m+i+3}(\a+1)_{m+i+3}}\frac{(1-x)^m}{m!}\frac{(2\a+5)_m}{(i+3)_m}\nn
&=&{n+\a \ch n}(-1)^{i+1}\frac{(i+2)!}{2^{i+3}}\sum_{m=0}^{\infty}
\frac{(-n)_{m+i+3}(n+2\a+1)_{m+i+3}}{(m+i+2)!(\a+1)_{m+i+3}}\frac{(2\a+5)_m}{m!}
\l(\frac{1-x}{2}\r)^m\n\eea
and therefore we obtain
\bea & &\sum_{i=2}^{\infty}c_i^*(x)D^{i+1}\SP\nn
&=&\frac{1}{2}{n+\a \ch n}\sum_{i=0}^{\infty}(-1)^{i+1}{\a+1 \ch i}\times{}\nn
& &\hspace{4cm}{}\times\sum_{m=0}^{\infty}
\frac{(-n)_{m+i+3}(n+2\a+1)_{m+i+3}}{(m+i+2)!(\a+1)_{m+i+3}}
\frac{(2\a+5)_m}{m!}\l(\frac{1-x}{2}\r)^m\nn
&=&-\frac{1}{2}{n+\a \ch n}\sum_{m=0}^{\infty}\frac{(-n)_{m+3}(n+2\a+1)_{m+3}}
{(m+2)!(\a+1)_{m+3}}\frac{(2\a+5)_m}{m!}\l(\frac{1-x}{2}\r)^m\times{}\nn
& &\hspace{5cm}{}\times\hyp{3}{2}{-n+m+3,-\a-1,n+2\a+m+4}{m+3,\a+m+4}{1}.\n\eea
By using the Saalsch\"{u}tz summation formula (\ref{Saal}) we have for $n\ge m+3$
\bea & &\hyp{3}{2}{-n+m+3,-\a-1,n+2\a+m+4}{m+3,\a+m+4}{1}\nn
&=&\frac{(\a+m+4)_{n-m-3}(-n-2\a-1)_{n-m-3}}{(m+3)_{n-m-3}(-n-\a)_{n-m-3}}\nn
&=&\frac{\G(n+\a+1)}{\G(\a+m+4)}\frac{(m+2)!}{\G(n)}\frac{\G(n+2\a+2)}
{\G(m+2\a+5)}\frac{\G(\a+m+4)}{\G(n+\a+1)}
=\frac{(m+2)!}{\G(n)}\frac{\G(n+2\a+2)}{\G(m+2\a+5)}\n\eea
which gives us by using (\ref{afgJ}), since $(-n)_{m+3}=0$ and the
$_3F_2(1)$ exists for $n<m+3$~:
\bea & &\sum_{i=2}^{\infty}c_i^*(x)D^{i+1}\SP\nn
&=&-\frac{1}{2}{n+\a \ch n}\sum_{m=0}^{\infty}\frac{(-n)_{m+3}(n+2\a+1)_{m+3}}
{(m+2)!(\a+1)_{m+3}}\frac{(2\a+5)_m}{m!}\times{}\nn
& &\hspace{8cm}{}\times\frac{(m+2)!}{\G(n)}\frac{\G(n+2\a+2)}
{\G(m+2\a+5)}\l(\frac{1-x}{2}\r)^m\nn
&=&-\frac{1}{2}{n+\a \ch n}\frac{\G(n+2\a+2)}{\G(n)\G(2\a+5)}\sum_{m=0}^{\infty}
\frac{(-n)_{m+3}(n+2\a+1)_{m+3}}{m!(\a+1)_{m+3}}\l(\frac{1-x}{2}\r)^m\nn
&=&\frac{4}{(2\a+3)(2\a+4)}{n+2\a+1 \ch n-1}\frac{d^3}{dx^3}\SP.\n\eea
This implies
\bea & &\sum_{i=2}^{\infty}c_i(x)D^{i+1}\SP=\frac{4}{(2\a+4)}{n+2\a+1 \ch n-1}
(1-x^2)\frac{d^3}{dx^3}\SP\nn
&=&\frac{4}{(2\a+4)}{n+2\a+1 \ch n-1}\l[2(\a+2)x\frac{d^2}{dx^2}\SP-
(n-1)(n+2\a+2)\frac{d}{dx}\SP\r].\n\eea

Now we have found
\bea & &\sum_{i=0}^{\infty}ic_i(x)D^i\SP+x\sum_{i=0}^{\infty}c_i(x)D^{i+1}\SP\nn
&=&4{n+2\a+1 \ch n-1}(1-x^2)\frac{d^2}{dx^2}\SP+4(2\a+3){n+2\a+2 \ch n-2}x\frac{d}{dx}\SP+{}\nn
& &\hspace{1cm}{}+\frac{4x}{(2\a+4)}{n+2\a+1 \ch n-1}\times{}\nn
& &\hspace{2cm}{}\times\l[2(\a+2)x\frac{d^2}{dx^2}\SP
-(n-1)(n+2\a+2)\frac{d}{dx}\SP\r]\nn
&=&4{n+2\a+1 \ch n-1}\frac{d^2}{dx^2}\SP\n\eea
which proves (\ref{second}) and therefore theorem 2.

\section{Differential equations for $P_n^{\a,\pm\frac{1}{2},0,N}(x)$}

In this section we will derive a differential equation for the polynomials
$\set{P_n^{\a,-\frac{1}{2},0,N}(x)}$ and another one for the polynomials
$\set{P_n^{\a,\frac{1}{2},0,N}(x)}$ for all $\a>-1$ and $N\ge 0$.
These differential equations can be obtained from our results by
applying the following quadratic transformations (see \cite{Koorn})~:
\be\la{tr1}\frac{P_{2n}^{\a,\a,M,M}(x)}{P_{2n}^{\a,\a,M,M}(1)}=
\frac{P_n^{\a,-\frac{1}{2},0,2M}(2x^2-1)}{P_n^{\a,-\frac{1}{2},0,2M}(1)}\ee
and
\be\la{tr2}\frac{P_{2n+1}^{\a,\a,M,M}(x)}{P_{2n+1}^{\a,\a,M,M}(1)}=
\frac{xP_n^{\a,\frac{1}{2},0,(4\a+6)M}(2x^2-1)}{P_n^{\a,\frac{1}{2},0,(4\a+6)M}(1)}.\ee
Note that (\ref{tr1}) and (\ref{tr2}) reduce to (\ref{even}) and (\ref{odd})
if $M=0$.

If we set $y(x):=f(2x^2-1)$ we can prove by induction that
$$y^{(2i)}(x)=\sum_{j=i}^{2i}\frac{(2i)!2^{3j-2i}}{(2j-2i)!(2i-j)!}x^{2j-2i}
f^{(j)}(2x^2-1),\;i=0,1,2,\ldots$$
and
$$y^{(2i+1)}(x)=\sum_{j=i+1}^{2i+1}\frac{(2i+1)!2^{3j-2i-1}}{(2j-2i-1)!(2i-j+1)!}
x^{2j-2i-1}f^{(j)}(2x^2-1),\;i=0,1,2,\ldots$$
or written in one formula
\be\la{y}y^{(i)}(x)=\sum_{j=\l[\frac{i+1}{2}\r]}^i\frac{i!2^{3j-i}}{(2j-i)!(i-j)!}
x^{2j-i}f^{(j)}(2x^2-1),\;i=0,1,2,\ldots.\ee

Note that (\ref{tr1}) and (\ref{y}) substituted in the differential equation
given in theorem 1 leads to a triviality, since
\bea\sum_{i=1}^{\infty}b_i(x)y^{(i)}(x)&=&\sum_{i=1}^{\infty}
\sum_{j=\l[\frac{i+1}{2}\r]}^i\frac{(-1)^i2^{3j-1}}{(2j-i)!(i-j)!}x^{2j}f^{(j)}(2x^2-1)\nn
&=&\sum_{j=1}^{\infty}2^{3j-1}x^{2j}f^{(j)}(2x^2-1)\sum_{i=j}^{2j}
\frac{(-1)^i}{(2j-i)!(i-j)!}=0,\n\eea
in view of
$$\sum_{i=j}^{2j}\frac{(-1)^i}{(2j-i)!(i-j)!}=\frac{(-1)^j}{j!}
\sum_{i=0}^j{j \ch i}(-1)^i=0,\;j=1,2,3,\ldots.$$

By using (\ref{tr1}) and (\ref{y}) we obtain from theorem 2 the
following equation~:
\bea & &M\sum_{i=0}^{\infty}c_i(x)\sum_{j=\l[\frac{i+1}{2}\r]}^i
\frac{i!2^{3j-i}}{(2j-i)!(i-j)!}x^{2j-i}f^{(j)}(2x^2-1)+{}\nn
& &\hspace{1cm}{}+(1-x^2)\l[16x^2f''(2x^2-1)+4f'(2x^2-1)\r]+{}\nn
& &\hspace{1cm}{}-8(\a+1)x^2f'(2x^2-1)+2n(2n+2\a+1)f(2x^2-1)=0\n\eea
with
$$c_0(x):=c_0(2n,\a,x)=4(2\a+3){2n+2\a+2 \ch 2n-2},\;\ndots,$$
satisfied by
$$f(x):=P_n^{\a,-\frac{1}{2},0,2M}(x).$$
Since
\bea & &\sum_{i=1}^{\infty}c_i(x)\sum_{j=\l[\frac{i+1}{2}\r]}^i
\frac{i!2^{3j-i}}{(2j-i)!(i-j)!}x^{2j-i}f^{(j)}(2x^2-1)\nn
&=&\sum_{j=1}^{\infty}2^{3j}f^{(j)}(2x^2-1)\sum_{i=j}^{2j}
\frac{i!2^{-i}}{(2j-i)!(i-j)!}x^{2j-i}c_i(x),\n\eea
we obtain
\bea\la{hulpdv1} & &M\sum_{j=0}^{\infty}d_j(x)f^{(j)}(2x^2-1)+
(1-x^2)\l[16x^2f''(2x^2-1)+4f'(2x^2-1)\r]+{}\nn
& &\hspace{1cm}{}-8(\a+1)x^2f'(2x^2-1)+2n(2n+2\a+1)f(2x^2-1)=0,\eea
where
$$\l\{\ba{l}\ds d_0(x)=4(2\a+3){2n+2\a+2 \ch 2n-2},\;\ndots\\ \\
\ds d_j(x)=\sum_{i=j}^{2j}\frac{i!2^{3j-i}}{(2j-i)!(i-j)!}x^{2j-i}c_i(x),
\;j=1,2,3,\ldots.\ea\r.$$
By using (\ref{cstar}) and (\ref{cgran}) we easily see that $d_j(x)$ is an
even polynomial with degree$\l[d_j(x)\r]\le 2j$ for each $j=1,2,3,\ldots$. Now we set $2x^2-1=t$ and $N=2M$
in (\ref{hulpdv1}) to find a differential equation of the form
\bea & &N\sum_{j=0}^{\infty}d_j^*(t)y^{(j)}(t)+(1-t^2)y''(t)+{}\nn
& &\hspace{1cm}{}-\frac{1}{2}\l[(2\a+1)+(2\a+3)t\r]y'(t)+\frac{1}{2}n(2n+2\a+1)y(t)=0,\n\eea
for the polynomials $\set{P_n^{\a,-\frac{1}{2},0,N}(t)}$, where
$$d_0^*(t)=\frac{1}{2}(2\a+3){2n+2\a+2 \ch 2n-2},\;\ndots$$
and
$$d_j^*(t)=\frac{1}{8}d_j\l(\sqrt{\frac{1+t}{2}}\r),\;j=1,2,3,\ldots.$$
We remark that $d_j^*(t)$ is a polynomial in $t$ with degree$\l[d_j^*(t)\r]\le j$ for
every $j=1,2,3,\ldots$. By using (\ref{cstar}) and (\ref{cgran}) we see that
\bea d_j^*(-1)&=&\frac{1}{8}d_j(0)=\frac{(2j)!2^{j-3}}{j!}c_{2j}(0)\nn
&=&(2\a+3)(-1)^{j+1}{\a+1 \ch j-1}\frac{2^{j-1}}{j!}\ne 0,\;j=1,2,3,\ldots\n\eea
if $\a$ is not a nonnegative integer. This implies that the order of the
differential equation is infinite in that case if $N>0$.
For nonnegative integer values of $\a$ we have $d_j(x)=0$ for $j>2\a+4$ and
$$d_{2\a+4}(x)=2^{4\a+8}x^{2\a+4}c_{2\a+4}(x)\ne 0$$
since degree$\l[c_{2\a+4}(x)\r]=2\a+4$.
This implies that the order of this differential equation equals $2\a+4$
if $\a$ is a nonnegative integer and $N>0$.

And if we set $y(x):=xf(2x^2-1)$ we find by using Leibniz' rule
\bea y^{(2i)}(x)&=&x\sum_{j=i}^{2i}\frac{(2i)!2^{3j-2i}}{(2j-2i)!(2i-j)!}
x^{2j-2i}f^{(j)}(2x^2-1)+{}\nn
& &\hspace{1cm}{}+2i\sum_{j=i}^{2i-1}\frac{(2i-1)!2^{3j-2i+1}}
{(2j-2i+1)!(2i-1-j)!}x^{2j-2i+1}f^{(j)}(2x^2-1)\nn
&=&\sum_{j=i}^{2i}\frac{(2i+1)!2^{3j-2i}}{(2j-2i+1)!(2i-j)!}x^{2j-2i+1}
f^{(j)}(2x^2-1),\;i=0,1,2,\ldots\n\eea
and
\bea y^{(2i+1)}(x)&=&x\sum_{j=i+1}^{2i+1}\frac{(2i+1)!2^{3j-2i-1}}
{(2j-2i-1)!(2i-j+1)!}x^{2j-2i-1}f^{(j)}(2x^2-1)+{}\nn
& &\hspace{1cm}{}+(2i+1)\sum_{j=i}^{2i}\frac{(2i)!2^{3j-2i}}{(2j-2i)!(2i-j)!}
x^{2j-2i}f^{(j)}(2x^2-1)\nn
&=&\sum_{j=i}^{2i+1}\frac{(2i+2)!2^{3j-2i-1}}{(2j-2i)!(2i-j+1)!}x^{2j-2i}
f^{(j)}(2x^2-1),\;i=0,1,2,\ldots\n\eea
or written in one formula
\be\la{yy}y^{(i)}(x)=\sum_{j=\l[\frac{i}{2}\r]}^i\frac{(i+1)!2^{3j-i}}
{(2j-i+1)!(i-j)!}x^{2j-i+1}f^{(j)}(2x^2-1),\;i=0,1,2,\ldots.\ee
If we substitute (\ref{tr2}) and (\ref{yy}) in the differential equation in
theorem 1 we find a triviality, since
\bea y(x)+\sum_{i=1}^{\infty}b_i(x)y^{(i)}(x)&=&xf(2x^2-1)+\sum_{i=1}^{\infty}
\sum_{j=\l[\frac{i}{2}\r]}^i\frac{(i+1)(-1)^i2^{3j-1}}{(2j-i+1)!(i-j)!}
x^{2j+1}f^{(j)}(2x^2-1)\nn
&=&\sum_{j=1}^{\infty}2^{3j-1}x^{2j+1}f^{(j)}(2x^2-1)
\sum_{i=j}^{2j+1}\frac{(i+1)(-1)^i}{(2j-i+1)!(i-j)!}=0,\n\eea
in view of
\bea\sum_{i=j}^{2j+1}\frac{(i+1)(-1)^i}{(2j-i+1)!(i-j)!}&=&\sum_{i=0}^{j+1}
\frac{(i+j+1)(-1)^{i+j}}{(j-i+1)!i!}\nn
&=&\frac{(-1)^j}{j!}\l[\sum_{i=0}^{j+1}{j+1 \ch i}(-1)^i-\sum_{i=0}^j{j \ch i}(-1)^i\r]\nn
&=&0,\;j=1,2,3,\ldots.\n\eea

By using (\ref{tr2}) and (\ref{yy}) we obtain from theorem 2 the
following equation~:
\bea & &M\sum_{i=0}^{\infty}c_i(x)\sum_{j=\l[\frac{i}{2}\r]}^i
\frac{(i+1)!2^{3j-i}}{(2j-i+1)!(i-j)!}x^{2j-i+1}f^{(j)}(2x^2-1)+{}\nn
& &\hspace{1cm}{}+(1-x^2)\l[16x^3f''(2x^2-1)+12xf'(2x^2-1)\r]+{}\nn
& &\hspace{1cm}{}-2(\a+1)x\l[4x^2f'(2x^2-1)+f(2x^2-1)\r]
+2(2n+1)(n+\a+1)xf(2x^2-1)=0\n\eea
with
$$c_0(x):=c_0(2n+1,\a,x)=4(2\a+3){2n+2\a+3 \ch 2n-1},\;\ndots,$$
satisfied by
$$f(x):=P_n^{\a,\frac{1}{2},0,(4\a+6)M}(x).$$
Since
\bea & &\sum_{i=0}^{\infty}c_i(x)\sum_{j=\l[\frac{i}{2}\r]}^i
\frac{(i+1)!2^{3j-i}}{(2j-i+1)!(i-j)!}x^{2j-i+1}f^{(j)}(2x^2-1)\nn
&=&\sum_{j=0}^{\infty}2^{3j}f^{(j)}(2x^2-1)\sum_{i=j}^{2j+1}
\frac{(i+1)!2^{-i}}{(2j-i+1)!(i-j)!}x^{2j-i+1}c_i(x),\n\eea
we obtain, after division by $x$
\bea\la{hulpdv2} & &M\sum_{j=0}^{\infty}e_j(x)f^{(j)}(2x^2-1)+
(1-x^2)\l[16x^2f''(2x^2-1)+12f'(2x^2-1)\r]+{}\nn
& &\hspace{1cm}{}-2(\a+1)\l[4x^2f'(2x^2-1)+f(2x^2-1)\r]+{}\nn
& &\hspace{1cm}{}+2(2n+1)(n+\a+1)f(2x^2-1)=0,\eea
where
$$e_j(x)=\sum_{i=j}^{2j+1}\frac{(i+1)!2^{3j-i}}{(2j-i+1)!(i-j)!}x^{2j-i}c_i(x),
\;j=0,1,2,\ldots.$$
By using (\ref{cstar}) and (\ref{cgran}) we easily see that $e_j(x)$ is an
even polynomial with degree$\l[e_j(x)\r]\le 2j$ for each $j=0,1,2,\ldots$. Now we set $2x^2-1=t$ and
$N=(4\a+6)M$ in (\ref{hulpdv2}) to find a differential equation of the form
\bea & &N\sum_{j=0}^{\infty}e_j^*(t)y^{(j)}(t)+(1-t^2)y''(t)+{}\nn
& &\hspace{1cm}{}-\frac{1}{2}\l[(2\a-1)+(2\a+5)t\r]y'(t)+\frac{1}{2}n(2n+2\a+3)y(t)=0,\n\eea
for the polynomials $\set{P_n^{\a,\frac{1}{2},0,N}(t)}$, where
$$e_0^*(t)=\frac{1}{2}{2n+2\a+3 \ch 2n-1},\;\ndots$$
and
$$e_j^*(t)=\frac{1}{8(2\a+3)}e_j\l(\sqrt{\frac{1+t}{2}}\r),\;j=1,2,3,\ldots.$$
Note that $e_j^*(t)$ is a polynomial in $t$ with degree$\l[e_j^*(t)\r]\le j$ for
every $j=1,2,3,\ldots$. By using (\ref{cstar}) and(\ref{cgran}) we see that
\bea e_j^*(-1)&=&\frac{1}{8(2\a+3)}e_j(0)=(-1)^{j+1}\frac{2^{j-1}}{j!}
\l[2{\a \ch j-1}(\a+1)+(2j+1){\a+1 \ch j-1}\r]\nn
&=&(2\a+5)(-1)^{j+1}{\a+1 \ch j-1}\frac{2^{j-1}}{j!}\ne 0,
\;j=1,2,3,\ldots\n\eea
if $\a$ is not a nonnegative integer. This implies that the order of the
differential equation is infinite in that case if $N>0$.
For nonnegative integer values of $\a$ we have $e_j(x)=0$ for $j>2\a+4$ and
$$e_{2\a+4}(x)=2^{4\a+8}x^{2\a+4}c_{2\a+4}(x)\ne 0.$$
This implies that the order of this differential equation equals $2\a+4$
for nonnegative integer values of $\a$ and $N>0$.

\vspace{5mm}

{\bf Acknowledgements.}
I wish to thank Francisco Marcell\'{a}n and his colleagues for their
invitation and hospitality during my visit to Madrid in June 1991.
During this stay in Madrid the proof of theorem 2 was completed.

Many thanks are due to my uncle Jan Koekoek for carefully reviewing earlier
versions of this manuscript and for his useful ideas and advices concerning
the last section of this paper.

\end{document}